\documentstyle[12pt]{article}

\newlength{\abstractwidth}
\renewcommand{\thefootnote}{\fnsymbol{footnote}}
\renewcommand{\thanks}[1]{\footnote{#1}} 
\newcommand{\starttext}{ \setcounter{footnote}{0}
\renewcommand{\thefootnote}{\arabic{footnote}}}

\newcommand{\be}{\begin{equation}}
\newcommand{\bea}{\begin{eqnarray}}
\newcommand{\eea}{\end{eqnarray}} \newcommand{\ee}{\end{equation}}

\def\ba{\begin{eqnarray}}
\def\ea{\end{eqnarray}}

\def\D{{\cal D}}

\def\K{{\cal K}}
\def\cM{{\cal M}}

\def\cX{{\cal X}}
\def\cL{{\cal L}}
\def\r{\rho}
\def\Aut{{\rm Aut}\,}

\def\ra{\rightarrow}

\def\o{\omega}
\def\Re{{\rm Re}}

\def\det{{\rm det}}

\def\log{\,{\rm log}\,}

\def\o{\omega}

\def\al{\alpha}
\def\b{\beta}
\def\g{\gamma}
\def\d{\delta}
\def\e{\varepsilon}

\def\l{\lambda}
\def\m{\mu}

\def\o{\omega}

\def\r{\rho}
\def\si{\sigma}
\def\t{\theta}

\def\D{\Delta}
\def\O{\Omega}

\def\na{\nabla}

\def\ti{\tilde}

\def\R{{\bf R}}
\def\C{{\bf C}}
\def\P{{\bf P}}

\def\i{\infty}

\def\p{\prod}
\def\s{\sum}

\def\ddb{{\partial\bar\partial}}

\def\ra{\rightarrow}

\def\D{\Delta}

\def\cM{{\cal M}}

\def\na{{\nabla}}

\def\K{{K\"ahler\ }}

 \def\v{\vskip .1in}

\def\[{{\bf [}}
\def\]{{\bf ]}}

\def\pl{\partial}

\begin{document}
\starttext \baselineskip=15pt \setcounter{footnote}{0}
\newtheorem{theorem}{Theorem}
\newtheorem{lemma}{Lemma}
\newtheorem{definition}{Definition}

\begin{center}
{\Large \bf THE DIRICHLET PROBLEM FOR DEGENERATE COMPLEX
MONGE-AMPERE EQUATIONS \footnote{Work supported in part by
DMS-07-57372 and DMS-05-14003.}}
\bigskip\bigskip

{\large  D.H. Phong and Jacob Sturm} \\

\end{center}

\medskip

\begin{abstract}

The Dirichlet problem for a Monge-Amp\`ere equation corresponding
to a non-negative, possible degenerate cohomology class on a K\"ahler
manifold with boundary is studied. $C^{1,\alpha}$ estimates away
from a divisor are obtained, by combining techniques of Blocki, Tsuji, Yau,
and pluripotential theory. In particular, $C^{1,\alpha}$ geodesic rays
in the space of K\"ahler potentials are constructed for each test
configuration.

\end{abstract}

\section{Introduction}
\setcounter{equation}{0}

This paper is motivated by two closely related problems: on one hand,
the Dirichlet problem for degenerate complex Monge-Amp\`ere equations on compact
K\"ahler manifolds with boundary,
and on the other hand, the existence
and regularity of geodesics in the space of K\"ahler potentials
of a given Chern class.

\smallskip

The complex Monge-Amp\`ere equation has been studied extensively over the years.
Three particularly influential works have been the work
by Yau \cite{Y78} on non-degenerate equations on compact K\"ahler manifolds without
boundary, by Caffarelli-Kohn-Nirenberg-Spruck \cite{CKNS}
on the Dirichlet problem on strongly pseudo-convex
domains in ${\bf C}^n$ (also for non-degenerate equations), and by Bedford-Taylor \cite{BT76,
BT82} on generalized solutions in the sense of pluripotential theory.
Since then, there has been considerable further progress, in particular
for $L^p$ right hand sides \cite{K98}, for existence and regularity
in the case
of quasi-projective manifolds (see e.g. \cite{CY80, CY86, K83, W07})
and for singular canonical metrics and degenerate equations
on K\"ahler manifolds
\cite{Ts, EGZ, TZ, DP, SW, ST, KT}. As for the Dirichlet problem, it has
been extended to more general domains in ${\bf C}^n$ \cite{G},
and to K\"ahler manifolds $M$ with boundary \cite{C00, reg, B09}
for equations of the form
\bea
\label{MA}
(\Omega+{i\over 2}\ddb\Phi)^m
=
F\,\Omega^m,
\quad \Omega+{i\over 2}\ddb\Phi\geq 0,
\qquad
\Phi_{\vert_{\pl M}}=\varphi,
\eea
when $F\geq 0$ and the given cohomology class  $[\Omega]$
is strictly positive definite (and hence is a K\"ahler form).
However, the important case of $[\Omega]$ non-negative, but
possibly degenerate on a K\"ahler manifold with boundary, has not been fully elucidated.

\smallskip
The degenerate complex Monge-Amp\`ere equation plays an important role in the
problem of constant scalar curvature metrics, which is a fundamental
problem in K\"ahler geometry.
Let $L\to X$ be a positive holomorphic line bundle over a compact complex manifold~$X$.
The positivity of $L$ means that there exists a metric $h_0$ on $L$ with
positive definite curvature $\o_0\equiv -{i\over 2}\ddb\log\,h$.
By a well-known conjecture of Yau \cite{Y93}, the existence of such a constant scalar curvature
metric $\o\in c_1(L)$
on a compact K\"ahler manifold $X$
with polarization $L\to X$, should be equivalent to the stability of
$L$ in the sense of geometric invariant theory. Suitable notions of K-stability
have been proposed since by Tian \cite{T97} and Donaldson \cite{D02}.
But of particular interest to us is yet another notion of stability introduced
by Donaldson \cite{D99}, where the one-parameter subgroups
of geometric invariant theory are replaced by geodesic rays in the
space ${\cal K}$ of K\"ahler potentials
\bea
\label{calK}
{\cal K}=\{\varphi\in C^\infty(X);\quad\o_0+{i\over 2}\ddb\varphi>0\}
\eea
equipped with the natural Weil-Petersson metric.
The conjecture/question of Donaldson \cite{D99} states that the absence of a constant scalar
curvature meteric should imply the existence of an infinite geodesic ray along which the
Mabuchi K-energy is decreasing.
\v

 If $(-T,0]\ni t\to
\varphi(\cdot,t)$ is a path in ${\cal K}$, Donaldson \cite{D99}, Semmes \cite{S}
and Mabuchi \cite{M} recognized that the geodesic equation
for $\phi$ is equivalent to the following completely degenerate Monge-Amp\`ere
equation
\bea
\label{HMA}
(\Omega_0+{i\over 2}\ddb\Phi)^m=0\ \ {\rm on}\ \ M,
\qquad
\Phi_{|w|=1}=\varphi(z,0)
\eea
where $M=X\times \{e^{-T}<|w|\leq 1\}$, $\Omega_0$ is the form $\o_0$ pulled back
to $M$, and
\bea
\Phi(z,w)\equiv\varphi(z,\log |w|),
\qquad (z,w)\in M.
\eea
We refer to $\varphi(\cdot,t)$ (equivalently $\Phi(\cdot,w)$)
as a {\it geodesic segment} when $T$ is finite, and as
a {\it geodesic ray} when $T$ is infinite. The Monge-Amp\`ere formulation
provides a notion of  ``generalized geodesics'': these are
rotation invariant, $\Omega_0$-plurisubharmonic,
generalized solutions $\Phi(z,w)$
of the equation (\ref{HMA}) in the sense of pluripotential theory.
Henceforth, we consider only generalized geodesics, and refer to them just as
geodesics for simplicity.

\smallskip

Using the Monge-Amp\`ere equation, geodesic segments in ${\cal K}$
have been constructed by Chen \cite{C00} using the method of a
priori estimates, and by \cite{PS06} using pluripotential theory. A
detailed analysis of the method of a priori estimates has been
provided very recently by Blocki \cite{B09}. A partial regularity
theory for the equation (\ref{HMA}) has been provided by \cite{CT}.
In this paper we are motivated by the stability condition of
Donaldson \cite{D99} and shall therefore focus on the regularity
problem for geodesic rays.

\smallskip

From the point of view of the Monge-Amp\`ere equation,
the difference between the segment and the ray cases resides
in the behavior of the ray near $w=0$, which is not prescribed by a Dirichlet
condition. For geodesic rays associated to a test configuration,
two possible approaches have been proposed in \cite{PS07} and \cite{reg}:
the geodesic ray can be modeled on approximating Bergman geodesics,
or it can be obtained as a solution of a Monge-Amp\`ere equation
on the compactification of the total space of the test configuration,
in a given cohomology class $[\Omega_0]$. We shall discuss this in greater detail
below. But here we note that \cite{reg} treated only the case of $[\Omega_0]$
strictly positive definite. The case where $[\Omega_0]$ is only assumed non-negative
is the topic of the present paper.

\smallskip
Our first results concern the general problem
of a priori estimates for the Dirichlet problem for the Monge-Amp\`ere equation
in the case of possibly degenerate cohomology class $[\Omega_0]$. Our precise set-up
is the following.

\smallskip

Let $M$ be a compact complex manifold of complex dimension $m$,
with $C^\infty$ boundary $\pl M$. Let $\Omega_0$ be a smooth, closed, non-negative
$(1,1)$-form on $ M$,
\bea
\label{1.5}
\Omega_0 \geq 0,
\eea
but not necessarily strictly positive. Assume the following key condition:
there exists an effective divisor
$E$ on $M$, supported away from the boundary $\pl M$, and an $\e>0$ such that
\bea
\label{E}
\Omega_0-\e [E]>0
\eea
where $[E]$ is the integration current on $E$. We fix such an $\e>0$ once and for all.
If $O(E)$ is the holomorphic line bundle defined by $E$,
and $\sigma(z)$ is the canonical section of $O(E)$ which vanishes on $E$,
condition (\ref{E}) is equivalent to the existence of a metric $H(z)$ on $O(E)$ satisfying the
condition
\bea
\label{MA0}
\Omega_0^{\e}\equiv \Omega_0+{i\over 2}\ddb \log\, H(z)^\e>0.
\eea
Note that the form $\Omega_0^{\e}$ is smooth, so that it defines a K\"ahler form on $M$.

\begin{theorem}
\label{degenerateapriori}
{\rm (a priori estimates)} Let $M$, $\Omega_0$, $E$ be defined as above,
and satisfying the conditions (\ref{1.5}) and (\ref{E}).
Let $\Phi\in C^2(M)\cap C^4(M\setminus \pl M)$ be a solution of the Dirichlet problem
\bea
\label{degenerate}
(\Omega_0+{i\over 2}\ddb\Phi)^m= F\,(\Omega_0^{\e})^m,
\quad
\Omega_0+{i\over 2}\ddb\Phi\geq 0,
\qquad
\Phi_{\vert_{\pl M}}=\varphi,
\eea
where $F> 0$ is a smooth function on $M$,
and $\varphi\in C^\infty(\pl M)$.  Let $-B$ be a lower
bound for the bisectional curvature of $\O_0^{\e}$, with $B\e>1$ and $b\geq 0$
an upper bound for the scalar curvature. Then we have,
with the covariant derivatives and norms written below
taken with respect to the metric
$\Omega_0^{\e}$:

\medskip
{\rm (i)} $C^1$ estimates:
\bea
\label{C1}
|\na\Phi(z)|\leq C_1\|\sigma(z)\|^{-\e A_1},
\qquad
z\in  M\setminus E,
\eea
Here $A_1$ is a constant which depends only on an upper bound for
$B$ and $\sup_M|\na F^{1\over m}|$.
The constant
$C_1$ depends on the dimension $m$ and upper bounds
for
$$
B,\ \sup_M|\na F^{1\over m}|,\ \sup_MF ,\  \sup_M|\Phi| ,\  \sup_M\|\sigma\|^{2\e} ,\
\sup_{\pl M}|\na\Phi|,\ \sup_{\pl M}{1\over \|\sigma\|^{2\e}}
$$

\medskip
{\rm (ii)} $C^2$ estimates:
\bea
\label{C2}
|\Delta\Phi(z)|\leq C_2\|\sigma(z)\|^{-\e A_2},
\qquad
z\in  M\setminus E.
\eea
Here the constant $A_2$ depends only on
upper bounds for $-\D\log F$,
$b$ and
$B$.
The constant $C_2$ depends on the dimension $m$ and upper bounds for
$$-\D\log F,\ b,\ B,\ {\rm osc}_M(\Phi),\ \sup_M\|\sigma\|^{2\e},\
\sup_{\pl M}(m+\D\phi)
$$

\end{theorem}
\v

We note that the constants $C_1, C_2, A_1, A_2$ do not depend on either
$\inf_M F$ or a lower bound for the bisectional curvature of $\Omega_0$,
an important fact which will be needed for the proofs of Theorem 2 and Theorem 3 below.
\v

The condition (\ref{E}) is a well-known type of condition in the literature.
It was introduced by Tsuji \cite{Ts},
who derived the type of $C^2$ estimate (ii) in Theorem \ref{degenerateapriori}
for the K\"ahler-Ricci flow on manifolds without boundary.
Recent related advances are in \cite{EGZ, TZ, DP, SW}. For manifolds with boundary,
$C^1$ estimates pose new difficulties. For strictly positive cohomology classes
$[\Omega_0]$, the $C^1$ estimates have been obtained by Chen \cite{C00} and, more recently,
by a completely different and more explicit method in the important advance of
Blocki \cite{B08}. For the case of degenerate
$[\Omega_0]$ considered here, we adapt the method of Blocki, with a key modification which shows
that the maximum principle can be applied at a definite distance from the divisor $E$.
We use the a priori estimates in Theorem \ref{degenerateapriori},
to  establish the following theorem:

\begin{theorem}
\label{degenerateexistence}
Let $M$, $\Omega_0$, $E$ be as above, and satisfy the conditions (\ref{1.5}), (\ref{E}).
Consider the Dirichlet problem for the totally degenerate
Monge-Amp\`ere equation
\bea
\label{MA00}
(\Omega_0+{i\over 2}\ddb\Phi)^m=0,
\quad
\Omega_0+{i\over 2}\ddb\Phi\geq 0
\qquad
\Phi_{\vert_{\pl M}}=0.
\eea
Then the Dirichlet problem
(\ref{MA00}) admits a bounded generalized solution $\Phi$
on the whole of $M$ in the sense of pluripotential theory.
Furthermore, for any $0<\al<1$, $\Phi$ is in $C^{\al}(M\setminus E)$.
If $M$ has a locally flat boundary, in the sense that near
each point on $\pl M$, there exists local holomorphic coordinates $(z_1,\cdots,z_m)$
so that $\pl M$ is given by $\Re\, z_m=0$,
then $\Phi\in C^{1,\al}(M\backslash E)$.
\end{theorem}

\medskip

We return now to the problem of constructing geodesic rays in the space
${\cal K}$ of K\"ahler potentials in a fixed Chern class $c_1(L)$,
where $L\to X$ is a positive line bundle (c.f. (\ref{calK})).
This problem has received
considerable attention recently,
and geodesic rays have been
constructed in \cite{AT, C00, PS07, reg, C08, CT07} by several different methods.

\smallskip

In particular, in \cite{PS07, reg}, it was shown how geodesic rays
can be constructed in two different ways starting from a test configuration for
$L\to X$. The open question is how to determine their asymptotic behavior as $t\to-\infty$.
More precisely, a test configuration
$\rho$ for $L\to X$ is an
endomorphism
\bea
\rho:{\bf C}^{\times}\to
{\rm Aut}({\cal L}\to{\cal X}\to {\bf C})
\eea
with the fibers of ${\cal L}\to{\cal X}$ over the point $w=1$
being isomorphic to $L\to X$ (see \cite{D02}). For geodesics, we restrict the fibration
to the closed unit disk $\{w\in {\bf C};|w|\leq 1\}$, with the point $w=0$ corresponding
to $t=-\infty$. The fiber $X_0$ of ${\cal X}$ over $w=0$
is called the central fiber.
Let ${\cal X}^\times$ be the complement of $X_0$ in ${\cal X}$.
In \cite{PS07}, the associated geodesic ray is a solution
of the complex Monge-Amp\`ere equation
on ${\cal X}^\times$,
obtained as a limit of Bergman geodesics,
and it was difficult to determine the regularity of the limiting ray.
Only for toric varieties has this problem been overcome very recently by
Song and Zelditch \cite{SZ06, SZ08}, using a refined semi-classical analysis
of norming constants.
In \cite{reg},
the associated geodesic ray was obtained from a priori estimates for a Monge-Amp\`ere equation
on the full space ${\cal X}$ (more precisely, the full space $\tilde{\cal X}$ of
an equivariant resolution,
\bea
p:\tilde{\cal X}\to {\cal X}\to {\bf C}.
\eea
of ${\cal X}$).
The ray is then of class $C^{1,\al}$ on $\tilde{\cal X}$
for any $0<\al<1$.
But the rays constructed in \cite{reg} are defined by a
strictly positive cohomology class on $\tilde{\cal X}$ (which is a kind of Dirichlet data
at the origin for the punctured disk), while certain more intrinsic cohomology classes
exist which are just non-negative. An important example is the class
of the pull-back $p^*({\cal L})$.
For such classes, a major question is whether the Dirichlet problem for the completely degenerate
Monge-Amp\`ere equation can be solved, and whether the solution is $C^{1,\al}$ away
from the central fiber.

\smallskip
Our next result is that, with Theorems \ref{degenerateapriori} and \ref{degenerateexistence},
this question
can be answered in the affirmative:

\smallskip

\begin{theorem}
\label{geodesicray}
Let $L\to X$ be a positive line bundle over a compact complex manifold $X$,
and let
$\rho:{\bf C}^{\times}\to
{\rm Aut}({\cal L}\to{\cal X}\to {\bf C})$ be a test configuration.
Fix a metric $h_0$ on $L$ with positive curvature $\o_0$.
Let $p:\tilde{\cal X}\to {\cal X}\to{\bf C}$ be an equivariant resolution
of ${\cal X}$, and let $\tilde{\cal X}^\times=p^{-1}({\cal X}^\times)$
be the complement of the central fiber in $\tilde{\cal X}$. Let $p^*{\cal L}$
be the pull-back to $\tilde{\cal X}$ of the line bundle ${\cal L}$.
Let $H_0$ be a metric on $p^*({\cal L}^k)$ with non-negative curvature
\bea
\Omega_0\equiv -{i\over 2}\ddb \,log\,H_0.
\eea
which satisfies the three properties enumerated in Lemma \ref{Omegazero}
below. Then there exists a generalized geodesic ray starting from $h_0$, i.e.,
there exists a bounded, rotation invariant,
$\Omega_0$-pluri\-subharmonic function $\Phi$ on $\tilde{\cal X}$ which is
a solution of the following Dirichlet problem,
\bea
(\Omega_0+{i\over 2}\ddb\Phi)^{n+1}=0\ \ {\rm on}\ \tilde{\cal X},
\qquad
\Phi_{\vert_{\pl\tilde{\cal X}}}=0,
\eea
in the sense of pluripotential theory. Furthermore, for any $0<\al<1$,
the function $\Phi$ is of class $C^{1,\al}$ in any compact subset of
$\tilde{\cal X}^{\times}$.
\end{theorem}

\section{Identities of Blocki, Yau, and Aubin}
\setcounter{equation}{0}

It is convenient to group together in this section
all the identities needed later for the derivation
of the a priori $C^1$ and $C^2$ estimates.
Thus consider the Monge-Amp\`ere equation
\bea
\label{MA}
(\Omega+{i\over 2}\ddb\Phi)^m\ =\ F\,\Omega^m
\eea
on an open manifold $M$, where $\Omega\equiv {i\over 2}g_{\bar kj}dz^j\wedge d\bar z^k$ is a
K\"ahler form (and in particular positive definite).
Let $-B$ be a lower bound for the bisectional
curvature of the K\"ahler
form $\Omega$, with $B\geq 0$. Thus
$R_{\bar j i \bar l k}a^i\overline{a^j}b^k\overline{b^l}\geq -B|a|^2|b|^2 $ for all
vectors $a^i, b^k$.
All covariant derivatives $\na$ in the identities below
are taken with respect to the K\"ahler form $\Omega$.
The form $\Omega'=\Omega+{i\over 2}\ddb\Phi$
is assumed to be strictly positive. The corresponding metric is denoted by
$g_{\bar kj}'$, and the corresponding covariant derivatives by $\na'$.
Unless otherwise indicated, norms of tensors are taken with respect to
$g_{\bar kj}$. But there are situations when $(g')_{\bar kj}$ or even a mixture
of $g_{\bar kj}$ and $(g')_{\bar kj}$ are also used, in which case
we shall use sub-indices to indicate which are the relevant metrics. For
example, we would write
\bea
|\na\na\Phi|_{\O\O'}^2=(g')^{j\bar k}g^{m\bar p}\na_j\na_m\Phi\na_{\bar k}\na_{\bar p}\Phi,
\eea
and so on.

\subsection{The Blocki identity}

The $C^1$ identity is due to Blocki \cite{B08}.
For a closer analogy with the identities of Yau, Aubin, and Calabi
given below, we formulate it in terms of covariant derivatives and
in terms of the endomorphism
$h=h^j{}_k$ defined by
\bea
\label{h}
h^j{}_k=g^{j\bar p}g_{\bar pk}'.
\eea
Let $\g:{\bf R}\to {\bf R}$ be a smooth real valued function,
and let $\Phi$ satisfy the Monge-Amp\`ere equation (\ref{MA}).
Set
\bea
\label{alpha}
\beta=|\na \Phi|^2_\O,
\qquad
\al= \log\,\beta-\g(\Phi).
\eea
Then at any interior critical point $p$ of $\al$, the following inequality \cite{B08} holds:
\bea
\label{blocki}
\Delta' \al
&\geq&
{1\over \beta}|\bar\na\na\Phi|_{\O\O'}^2
+
(\g'(\Phi)-B-{F_1\over \beta^{1\over 2}})\,{\rm Tr}\,h^{-1}
+
(-\g''(\Phi)+2{\g'(\Phi)\over \beta})|\na\Phi|_{\O'}^2
\nonumber\\
&&
\quad
-(m+2)\g'(\Phi)-{2\over\beta}
\eea
where the constant $F_1$ is defined by
\bea
F_1= 2\,{\rm sup}_M|\na (F^{1\over m})|.
\eea

For the convenience of the reader, we give the derivation of the inequality (\ref{blocki}),
following \cite{B09}. First, a direct calculation gives
\bea
\label{directcalculation}
\Delta'\al
&=&
{1\over \b}\bigg\{|\na\na\Phi|_{\O\O'}^2+|\bar\na\na\Phi|_{\O\O'}^2+
2\Re (\pl^p(\log F)\pl_p\Phi)+(g')^{k\bar l}\pl_p\Phi R^p{}_{k\bar l}{}^{\bar m}\pl_{\bar
m}\Phi\bigg\}
\nonumber\\
&&
-{1\over \b^2}|\na\b|_{\O'}^2
-\g''(\Phi)|\na\Phi|_{\O'}^2
-\g'(\Phi)\Delta'\Phi,
\eea
Here we make use of the identity $(g')^{p\bar q}\na_m g_{\bar qp}'=\pl_m\log F$, which is obtained
by differentiating the Monge-Amp\`ere equation (\ref{MA}). It is readily seen that
\bea
2|\pl^p(\log F)\pl_p\Phi|&\leq & 2|\na (\log F)|\b^{1\over 2} \leq \
m{F_1\over F^{1\over m}}\beta^{1\over 2} \
\leq \ ({\rm Tr}\,h^{-1})\, F_1\b^{1\over 2}
\nonumber\\
(g')^{k\bar l}\pl_p\Phi R^p{}_{k\bar l}{}^{\bar m}\pl_{\bar m}\Phi
&\geq&
-B\beta {\rm Tr}\,h^{-1}
\nonumber\\
-\Delta'\Phi&=&-m+{\rm Tr}\,h^{-1},
\eea
so the only troublesome term is $-\beta^{-2}|\na\b|_{\O'}^2$ (if $\g'$, $-\g''$
are chosen to be non-negative, as will be the case later in applications).

\smallskip
The idea in \cite{B09} is to cancel this term partially with $|\na\na\Phi|_{\O\O'}^2$. First,
introduce the quantity
\bea
\Psi_p\equiv g^{j\bar k}\na_j\na_p\Phi\na_{\bar k}\Phi.
\eea
On one hand, we have, by the Cauchy-Schwarz inequalty,
\bea
\beta |\na\na\Phi|_{\O\O'}^2\geq (g')^{p\bar q}\Psi_p\overline{\Psi_q},
\eea
on the other hand, by a simple calculation
\bea
\na_p\b=\Psi_p+\na_j\Phi h^j{}_p-\na_p\Phi.
\eea
Altogether, we find
\bea
{1\over\b}|\na\na\Phi|_{\O\O'}^2
\geq |{\na\b\over\b}-{1\over \b}\na\Phi\, h+{1\over\b}\na\Phi|_{\O'}^2.
\eea

\smallskip
So far, all the calculations have been at an arbitrary point. Now, assume that we
are at an interior critical point of the function $\al$. Then
\bea
{\na\b\over\b}=\g'(\Phi)\na\Phi
\eea
and the preceding expression simplifies (again, with the assumption $\g'>0$)
\bea
|{\na\b\over\b}-{1\over \b}\na\Phi\, h+{1\over\b}\na\Phi|_{\O'}^2
&=&
|\g'(\Phi)\na\Phi-{1\over \b}\na\Phi\, h+{1\over\b}\na\Phi|_{\O'}^2
\nonumber\\
&\geq&
(\g'(\Phi)^2)|\na\Phi|_{\O'}^2
-
2\g'(\Phi)+2{\g'(\Phi)\over\beta}|\na\Phi|_{\O'}^2-{2\over \b}.
\eea
Substituting this inequality in (\ref{directcalculation}) gives the desired estimate.

\subsection{Yau and Aubin identities}

The identity of Yau and Aubin is the following, formulated in terms
of the endomorphism $h^j{}_k$ as in \cite{PSS},
\bea
\label{yau}
\Delta'\,\log\,{\rm Tr}\,h
&=&
{1\over {\rm Tr}\,h}(-R+\Delta F+(h^{-1})^p{}_m R^m{}_p{}^j{}_kh^k{}_j)
\nonumber\\
&&
\quad
+
\bigg\{
(g')^{p\bar q}{\rm Tr}(\na_p'h \,h^{-1}\na_{\bar q}'h)
-
{|\na'{\rm Tr}\,h|^2\over ({\rm Tr}\,h)^2}\bigg\}.
\eea
Again, we provide the derivation. Write
\bea
\Delta' {\rm Tr} h
&=&(g')^{p\bar q}\pl_{\bar q}\pl_p{\rm Tr} h
=
(g')^{p\bar q}{\rm Tr}(\na_{\bar q}'[(\na_p' h\,h^{-1})h])
\nonumber\\
&=&
(g')^{p\bar q}{\rm Tr}(\na_{\bar q}'(\na_p' h \,h^{-1})h)
+(g')^{p\bar q}{\rm Tr}(\na_p' h\,h^{-1}\na_{\bar q}'h).
\eea
Now the curvatures of $\O$ and $\O'$ are related by the standard
formula
$$
\na_{\bar q}'(\na_p'h\,h^{-1})=-(Rm)_{\bar qp}(\O')+
(Rm)_{\bar qp}(\O),
$$
where $(Rm)_{\bar qp}$
denotes the Riemann curvature tensor viewed as an endomorphism
of the tangent bundle. Substituting this identity in the previous one gives
\bea
\Delta'{\rm Tr} h
=
(g')^{p\bar q}{\rm Tr}(-Rm'_{\bar qp}\cdot h +Rm_{\bar qp}\cdot h)
+
(g')^{p\bar q}{\rm Tr}(\na_p' h\,h^{-1}\na_{\bar q}'h),
\eea
Thus
\bea
\Delta'\log{\rm Tr}\,h
&=&
{\Delta' {\rm Tr}\, h\over {\rm Tr}\, h}
-
{|\na'{\rm Tr}\, h|_{\O'}^2
\over({\rm Tr}\,h)^2}
\\
&=&
{(g')^{p\bar q}{\rm Tr}(-Rm_{\bar qp}'\cdot h+Rm_{\bar qp}\cdot h)\over {\rm Tr} h}
+
\big\{{(g')^{p\bar q}{\rm Tr}\,(\na_p'h \,h^{-1}\na_{\bar q}'h)\over {\rm Tr}\, h}
-
{|\na'{\rm Tr}\,h|_{\O'}^2\over({\rm Tr}\, h)^2}\big\}
\nonumber
\eea
Now $(g')^{p\bar q}=(h^{-1})^p{}_{m}g^{m\bar q}$, and the middle term on the right hand side
can be rewritten as
\bea
(g')^{p\bar q}{\rm Tr}(Rm_{\bar qp}\cdot h)
=(h^{-1})^p{}_m R^m{}_p{}^\al{}_\b h^\b{}_\al.
\eea
As for the first term on the right hand side, it can be recognized as
\bea
(g')^{p\bar q}{\rm Tr}(Rm'_{\bar qp}\cdot h)
=(Ric(\O'))^\al{}_\b h^\b{}_\al
=
(g')^{\al\bar\g}(R')_{\bar\g\b}h^\b{}_\al
=
g^{\b\bar\g}(R')_{\bar\g\b}.
\eea
The Monge-Amp\`ere equation implies that
\bea
(R')_{\bar\g\b}=R_{\bar\g\b}-\pl_{\bar\g}\pl_\b\log F,
\eea
and thus
\bea
g^{\b\bar\g}(R')_{\bar\g\b}=R-\Delta \log F.
\eea
Putting all together, we obtain the desired identity.

\medskip\noindent
Combining Yau's basic inequality
\bea
{(g')^{p\bar q}{\rm Tr}(\na_p'h \,h^{-1}\na_{\bar q}'h)\over {\rm Tr} \,h}
-
{|\na'{\rm Tr}\,h|^2\over ({\rm Tr}\,h)^2}\geq 0
\eea
with the simple identity used earlier
$\Delta'\Phi=m-{\rm Tr}\,h^{-1}$,
and the fact that $({\rm Tr}\,h)^{-1}\leq {\rm Tr}\,h^{-1}$,
we obtain
\bea\label{YA}
\Delta'(\,\log\,{\rm Tr}\,h-A\Phi)
\geq {A\over 2} {\rm Tr}\,h^{-1}-mA,
\eea
for all constants $A$ larger than a constant depending only on lower bounds for
$-R$ and $\Delta F$ and the bisectional curvature $R_{\bar jj\bar kk}$ of the metric $\Omega$.

\section{Proof of Theorem 1}
\setcounter{equation}{0}

We can now give the proof of Theorem \ref{degenerateapriori}. First, define $\Omega_0^{\e}$
by
\bea
\Omega_0^\e=\Omega_0+{i\over 2}\ddb \,\log\,H^\e
\eea
By hypothesis, the form $\Omega_0^\e$ is a K\"ahler form for some strictly positive
$\e$, which we fix from now on.
We use the method
of Tsuji \cite{Ts}
and rewrite the original equation (\ref{degenerate}) with $\Omega_0^\e$ as background,
\bea
\label{MA0e}
(\Omega_0^\e+{i\over 2}\ddb \Phi^\e)^m
= F \,(\Omega_0^\e)^m,
\eea
where we have defined $\Phi^\e$ by
\bea
\Phi^\e(z)=\Phi-\e\log \|\sigma(z)\|^2.
\eea
The function $\Phi^\e$ is smooth in $ M\setminus E$.

\subsection{Proof of the $C^1$ bound}

We now apply Blocki's identity (\ref{blocki}), with a suitable choice of function $\g(x)$.
For convenience, we refer to constants satisfying the dependence spelled out
in part (i) of Theorem~1 for the constant $C_1$ as ``admissible constants''.
Since $\Phi$ is continuous on $M$, the function $\Phi^\e(z)$
is bounded from below. Fix an admissible constant $C$
so that
\bea
\inf_M\,\Phi^\e\geq -C+1.
\eea
Let $-B$, with $B\geq 0$, be a lower bound for the bisectional curvature of the metric
$\Omega_0^\e$ and assume $B\e>1$.
Set as before $F_1^\e=2{\rm sup}|\na F^{1\over m}|$.
Define the function $\g(x)$ to be
\bea\label{gamma}
\g(x)
=
(B+F_1^\e+1)x-{1\over x+C}.
\eea
Note that, unlike in \cite{B08}, this choice of $\g(x)$ does not involve an {\it upper bound}
for $\Phi^\e$, which is not available in the present case. The function $\g(x)$
is smooth and bounded from below in the range $[-C+1,\infty)$, $\g'(x)>B+F_1^\e+1$
throughout this range, $\g''(x)=-2(x+C)^{-3}$, and $\g(x)\to +\infty$
as $x\to +\infty$.

\medskip

Consider the function $\al$ of (\ref{alpha}),
defined now for the equation (\ref{MA0e}),
that is, with $\Omega_0$ and $\Phi$ replaced now by $\Omega_0^\e$ and $\Phi^\e$
respectively,
\bea
\b=|\na\Phi^\e|^2,
\qquad
\al=\log\,\b-\g(\Phi^\e).
\eea
Since $|\Phi|$, $|\na\Phi|$ are all assumed to be continuous on $ M$,
the function $\al$ tends to $-\infty$ near the divisor $E$ (due to the requirement that
$B\e>1$).  Thus
it admits a maximum at a point $p$ which is either on $\pl M$, or in the interior of $M\setminus
E$. If $p$ is on $\pl M$, then the desired bound for $|\na\Phi|$ follows at once.

\medskip
Assume then that $p$ is an interior point of $M\setminus E$.
It suffices to show that there exists an admissible constant $C_3$, so
that
\bea
\label{alp}
\al(p) \leq C_3.
\eea
Indeed, this would imply that at an arbitrary point $z\in  M\setminus E$, we have
\bea
|\na\Phi^\e(z)|^2\leq C_3 \,e^{\g(\Phi^\e(z))}\ \leq C_3 e^{(B+F_1^\e+1)\Phi^\e(z)}
\leq C_4 \|\sigma(z)\|^{-2\e(B+F_1^\e+1)}.
\eea
Since $\na\Phi^\e=\na\Phi+\e\na\log\|\sigma\|^2$, we obtain the desired inequality
for $|\na\Phi|$,
with the exponent $A_1=2\e(B+F_1+1)$, since $B\e>1$.

\medskip
We turn to the proof of (\ref{alp}). Clearly, we can assume that $\beta>1$ at
the point $p$, since otherwise the boundedness from below of $\g$ would imply that $\al(p)$
satisfies the desired bound. With this assumption, Blocki's estimate
(\ref{blocki})
first simplifies to
\bea
\Delta' \al(p)
&\geq&
(\g'(\Phi^\e)-B-F_1)\,{\rm Tr}\,h^{-1}
-\g''(\Phi^\e)|\na\Phi^\e|_{\O'}^2
\nonumber\\
&&
\quad
-(m+2)\g'(\Phi^\e)-2\ ,
\eea
Here we have dropped the manifestly positive terms $|\bar\na\na\Phi^\e|^2$,
$\g'(\Phi^\e)\beta^{-1}|\na\Phi^\e|^2$.
Next, making use of the definition (\ref{gamma})  of $\g$,
the inequality simplifies further to
\bea
\label{C1maximum}
\Delta'\al(p)
&\geq&
{\rm Tr} \,h^{-1}+{2\over (\Phi^\e+C)^3}|\na\Phi^\e|_{\O'}^2
-(m+2)(B+F_1^\e+1)-(m+4).
\eea
Since $p$ is a maximum, $\Delta'\al(p)\leq 0$.
It follows that ${\rm Tr}\,h^{-1}(p)\leq C_5$. If $\lambda_i$
are the eigenvalues of $g_{\bar kj}'$ with respect to $g_{\bar kj}$ at the point $p$,
this implies that all
$\lambda_i^{-1}$ are bounded by admissible constants. On the other hand,
$\p_{j=1}^m\l_j=F$ so the $\l_j$ are also bounded above by admissible constants.
In particular, we have
\bea
|\na\Phi^\e(p)|_{\O}^2\leq C_6 |\na\Phi^\e(p)|_{\O'}^2
\eea
where $C_6$ is an admissible constant.

\medskip

The next  step is to show that $\|\sigma(p)\|^{-2\e}$
is bounded from above by an admissible constant. To see this,
we return to the inequality
(\ref{C1maximum}) and deduce that
\bea
|\na\Phi^\e(p)|^2\leq C_7\,(\Phi^\e(p)+C)^3
\leq C_7(C_8+\e\log{1\over \|\sigma(p)\|^2})^3.
\eea
Now, by the definition of admissible constants, we have
\bea
\al(z)\geq -C_9,
\qquad z\in \pl M
\eea
where $C_9$ is admissible. On the other hand,
$\g(x)\geq {1\over C_{10}}\,x-C_{10}$ with an admissible $C_{10}>0$, so we have
\bea
\g(\Phi^\e(p))
\geq {1\over C_{10}} \,\e\log{1\over\|\sigma(z)\|^2}-C_{11},
\eea
We thus obtain
\bea
-C_9\leq \al(p)
\leq
3\,\log (C_8+\e\log{1\over \|\sigma(p)\|^2})-{1\over C_{10}}
\,\e\log{1\over\|\sigma(z)\|^2}+C_{12}.
\eea
This shows that
\bea
\e\log{1\over\|\sigma(p)\|^2}\leq C_{13},
\eea
and hence $\al(p)\leq C_{13}$. This is the desired conclusion,
and the proof of the $C^1$ estimates is complete.

\subsection{Proof of the $\Delta\Phi$ bound}

For this section, the word ``admissible'' will now refer to
constants with the dependence spelled out for the constant $C_2$ in (ii).
We consider the expression ${\rm Tr}\,h-A\Phi^\e$.
Since
\bea
{\rm Tr}\,h&=&g^{j\bar k}(g')_{\bar kj}=m+\Delta \Phi^\e
=
m+\Delta\Phi-\e\Delta \log\|\sigma(z)\|^2
\nonumber\\
&=&
m+\Delta\Phi-\e\Delta\log\, H,
\eea
this expression is $C^2$ on $ M\backslash \pl M$ and extends to a $C^0$ function on $M$.
On the other hand, $\Phi^\e
=
\Phi-\e\log\|\sigma(z)\|^2\to +\infty$ as $z\to E$.
It follows that the expression ${\rm Tr}\,h-A\Phi^\e$ attains its maximum
either on the boundary $\pl M$ or an an interior point $p\in M\setminus E$.
It suffices to show that
\bea
(\log{\rm Tr}\,h-A\Phi^\e)(p)\leq C_{14}
\eea
for some admissible constant $C_{14}$.

\medskip

By the definition of admissible constants, this inequality holds if $p$ is on $\pl M$.
Thus assume that $p$ is in the interior of $M\setminus E$.
Since we have then $\Delta'(\log\,{\rm Tr}h-A\phi)(p)\leq 0$,
it follows from the Yau and Aubin
estimate (\ref{YA}) that
${\rm Tr}\, h^{-1}(p) \leq C_{15}$.
This implies that each eigenvalue $\lambda_i$ of $h$ is bounded from
below by an admissible constant. Since their product is $F$,
each of them is also bounded from above by an admissible constant.
For any $z\in M\setminus E$, we can then write
\bea
\log\,{\rm Tr}\,h(z)
\leq A (\Phi^\e(z)-\Phi^\e(p))
+
\log\,{\rm Tr}\, h(p)
\leq
A(\Phi^\e(z)-\Phi^\e(p))+C_{16}.
\eea
Thus
\bea
{\rm Tr}\,h(z)
&\leq& C_{17}\,e^{A(\Phi^\e(z)-\Phi^\e(p))}
\nonumber\\
&=&
C_{17}e^{A(\Phi(z)-\Phi(p))}\|\sigma(z)\|^{-2\e A}\|\sigma(p)\|^{2\e A}
\leq
C_{18}e^{A\,{\rm osc}(\Phi)}\|\sigma(z)\|^{-2\e A}.
\eea
The proof of the $C^2$ estimates is complete.

\section{Proof of Theorem 2}
\setcounter{equation}{0}

We wish to solve the Dirichlet problem
\bea
\label{MA00-1}
(\Omega_0+{i\over 2}\ddb\Phi)^m=0,
\qquad
\Omega_0+{i\over 2}\ddb\Phi\geq 0,
\qquad
\Phi_{\vert_{\pl M}}=0.
\eea
where the $(1,1)$-form $\Omega_0$ is smooth and non-negative, but may be degenerate.

\subsection{A viscosity approximation}

As in viscosity methods, we try and
obtain a solution to the above generate problem
as the limit of a subsequence of solutions
of the Dirichlet problem for non-degenerate equations.
More specifically, set
\bea
\label{Omegas}
\Omega_s\equiv \Omega_0+s{i\over 2}\ddb \log\,H^\e,
\qquad 0\leq s\leq 1.
\eea
Then $\Omega_s\to \Omega_0$ as $s\to 0$,
and for each $0<s\leq 1$, the form $\Omega_s$ is a K\"ahler form. In fact,
if $\Omega^\e\equiv \Omega_0+{i\over 2}\ddb\log\,H^\e$ is the K\"ahler form
constructed earlier,
we can write
\bea
\Omega_s=(1-s)\Omega_0+s\Omega^\e,
\eea
which  immediately implies
 that $\Omega_s$ is non-degenerate for $0<s\leq 1$.

\smallskip
Consider now for each sufficiently small $s>0$ the following Dirichlet problem
for the non-degenerate complex Monge-Amp\`ere equation,
\bea
\label{MA00-2}
(\Omega_s+{i\over 2}\ddb\Phi_s)^m=F_s\,\Omega_s^m,
\qquad
\Omega_s+{i\over 2}\ddb\Phi\geq 0,
\qquad
(\Phi_s){\vert_{\pl M}}=0,
\eea
where the right hand sides $F_s$ are smooth scalar functions on $\bar M$
which are defined as follows.
Let
\bea
\Omega_s^\e\equiv \Omega_s+{i\over 2}\ddb\,\log\,H^\e\ = \
\Omega^\e+s{i\over 2}\ddb\log H^\e
\eea
Then $\Omega_s^\e>0$ for $s>0$ sufficiently small. Then choose $F_s$
so that
\bea
\label{Fcondition1}
{\rm sup}_{0<s\leq 1}\|F_s\|_{C^0(\bar M)}\leq 1,\ \
\lim_{s\to 0}\,{\rm sup}_{0<s\leq 1}\|F_s\|_{C^0(\bar M)}\ = \ 0,
\eea
and such that
\bea\label{46}
F_s^{(\e)}=F_s\, {(\Omega_s)^m\over(\Omega_s^\e)^m}\
\eea
is a positive constant.

\subsection{Convergence to a generalized solution}

Since $\Omega_s$ is a K\"ahler form for $0<s<1$, and since,
by virtue of (\ref{Fcondition1}), the function
$\underline\Phi_s=0$ is a subsolution for the equation
(\ref{MA00-2}),
it follows from Theorem 1.3 of \cite{B09} that the Dirichlet problem
(\ref{MA00-2}) admits a strictly $\Omega_s$-plurisubharmonic solution
$\Phi_s$ which is smooth on $M$.

\medskip
{\it (a) Uniform $C^0$ estimates for $\Phi_s$}

\smallskip

We claim that there exists a constant $C$ independent of $s$, so
that
\bea
\label{C0}
\|\Phi_s\|_{C^0( M)}
\leq C.
\eea
First, by the $\Omega_s$-plurisubharmonicity condition,
we have
\bea
\Delta_{\Omega^\e}\Phi_s\geq -(\Omega^\e)^{j\bar k}(\Omega_s)_{\bar kj}
\geq -C
\eea
where $C$ is a constant independent of $s$. If $\hat\Phi$ is the
solution of the Dirichlet problem
\bea
\Delta_{\Omega^\e}\hat\Phi=-C,
\qquad \hat\Phi_{\vert_{\pl M}}=0,
\eea
it follows from the maximum principle for the Laplacian that
\bea
\Phi_s \leq \hat \Phi
\eea
for all $0<s\leq 1$. Next, recall that the function $0$ is a subsolution
of the Dirichlet problem (\ref{MA00-2}), in view of the fact
that the right hand side $F_s$ satisfies the condition (\ref{Fcondition1}).
It follows from the maximum principle for the Monge-Amp\`ere equation that
\bea
0\leq \Phi_s,
\qquad {\rm for\ all}\ 0<s\leq 1.
\eea
The estimate (\ref{C0}) follows.

\medskip
{\it (b) Uniform $C^1$  estimates for $\Phi_s$ at $\pl M$}.

\smallskip

Since $0\leq \Phi_s\leq \hat\Phi$, and $\hat\Phi$ vanishes at $\pl M$,
it follows that the absolute values of the partial derivatives of $\Phi_s$ are bounded
uniformly, where here, as before, all covariant derviatives and norms are
taken with respect to $\O_0^\e$.

\medskip
{\it (c) Uniform $C^1$  estimates on compact subsets of $M\setminus E$}\ .

\smallskip
We shall show the convergence
in $C^{\al}$
 over  $ M\setminus E$,
of a subsequence of $\Phi_s$ to
a bounded solution
the original equation (\ref{MA00-1}) in the sense of pluripotential theory. For this, we
shall apply
Theorem \ref{degenerateapriori} to obtain estimates for the gradient
of $\Phi_s$ which
is independent of $s$. We check the hypotheses of Theorem
\ref{degenerateapriori} in the present case: First note that
the equation (\ref{MA00-1}) can be readily put in the form of the
equations considered in Theorem \ref{degenerateapriori}
\bea
\label{MAFs}
(\Omega_s+{i\over 2}\ddb\Phi_s)^m
= F_s^{(\e)}(\Omega^\e)^m,
\eea
with the constant function $F_s^{(\e)}$ is defined by (\ref{46}).

\v
Recall that the functions $\Phi_s$
are $C^\infty$ in $\bar M$, with $C^0$ norms uniformly bounded in $s$,
by (\ref{C0}).
Applying Theorem \ref{degenerateapriori} we deduce the exitence of a
constants $C,A_1>0$, which are independent of $s$ and which satisfy
\be
\label{C1-C2}
|\na\Phi_s(z)|\leq  C\,\|\sigma(z)\|^{-2\e A_1}.
\ee

\smallskip This gives the desired uniform $C^1$ estimate.

\medskip
{\it (d) Uniform $C^{1,\al}$  estimates on compact subsets of
$M\setminus E$ in the case where $M$ has locally flat boundary.}

\v
Next we
derive the $C^{1,\al}$ estimate under the assumption of a locally flat boundary:
For non-degenerate Monge-Amp\`ere equations, a priori estimates at
the boundary depending on a lower bound for the right hand side
were obtained by Caffarelli et al. \cite{CKNS} and Guan \cite{G} by barrier arguments.
As noted by Blocki \cite{B09}, the same arguments can provide bounds
independent of a lower bound for the right hand side if the boundary
is locally flat. The same observation was used implicitly earlier by Chen \cite{C00}.
The following lemma contains a precise formulation
of the a priori estimate that we need:

\begin{lemma}\label{C2barrier}
 Let $M$ be a compact complex manifold of dimension $m$, with smooth locally flat boundary,
let $U$ be an
open neighborhood of  $\pl M$,  $\O$ a \K metric on $U$ and
let $F\in C^\i(U)$ a positive smooth function.  Suppose $\Phi\in C^\i(U)$ with
$\O'=\O+{i\over 2}\ddb\Phi>0$ and $(\O+{i\over 2}\ddb\Phi)^{m}=F\O^{m}$. Then there is a
constant
$C>0$, depending on $m,\O$, and upper bounds for $\sup_UF$ and $\sup_U|\na\log F|$, such that

\be \sup_{\pl M} (m+\D\Phi)\  \leq \
C\sup_{\pl M}(1+|\na\Phi|^2)\cdot\sup_{U}(1+|\na\Phi|^2)
\ee
where the norms and covariant derivatives are taken with respect to $\O$.
\end{lemma}
\smallskip
For the sake of completeness we provide a self-contained proof,
which is a straightforward adaptation of the arguments in
Guan \cite{G} and Chen \cite{C00}.
\v

{\it Proof of Lemma \ref{C2barrier}}: We write $\O={i\over 2}g_{\bar\b\al}dz^\al\wedge d\bar z^\b$
and $\O'\ ={i\over 2} g'_{\bar\b\al}dz^\al\wedge d\bar z^\b$.
Let $p\in \pl M$ be a point where
$\sup_{\pl M} (m+\D\Phi)$ is achieved, and choose   coordinates
$(z_1,...,z_{m})$ on $U$, centered at $p$ and $\d=\d_g>0$ so that
$g_{\bar\b\al}(0)=\d_{\bar\b\al}$
and so that
$$ {1\over 2}I\ \leq \ g(z) \ \leq
2I
$$
for all $z$ in the neighborhood
$U_\d=\{(z_1,...,z_{m}): \s |z_j|^2< \d, x=Re(z_{m})\geq 0\} $
\v

It suffices to show that
for $\al<m$ or $\b<m $, we have

\be\label{2.8} \left|
{\pl^2\Phi\over \pl z_\al\pl z_{\bar \b}}(0)
\right|\ \leq \ C(\sup_{\pl M} |\na\Phi| + 1)(\sup_{U} |\na\Phi| + 1)
\ee
To see this, first note we may assume $\al=m$ or $\b=m$, for
otherwise, the left side of (\ref{2.8}) vanishes identically since $\Phi$ vanishes
on the boundary. Next we observe that
$$ \det(\d_{\bar\b\al}+\pl_\al\pl_{\bar \b}\Phi)\ = \ F\cdot \det(g)
$$
so
$$
 1\ + \ {\pl^2\Phi\over \pl z_{m}\pl \bar z_{m}}(0)\ = \ F \ + \
\s_{\al=1}^{m-1}
{\pl^2\Phi\over \pl z_{\al}\pl \bar z_{m}}\cdot
{\pl^2\Phi\over \pl z_{m}\pl \bar z_{\al}}
$$
and thus Lemma \ref{C2barrier} follows from (\ref{2.8}).

\v
Now we prove (\ref{2.8}):
we assume $\al<m$ and $\b=m$.
Let $D=\pm {\pl\over \pl x_\al}$ or $\pm {\pl\over \pl y_\al}$
where $x_\al=Re(z_{\al}), y_\al = Im(z_{\al})$ and $\al<m$.
Let $x=Re(z_{m})$.
Then we must show
\be\label{cl} \left({\pl\over \pl x} D\Phi\right)(0)\ \leq \
C(\sup_{\pl M} |\na\Phi| + 1)(\sup_{U} |\na\Phi| + 1)
\ee
Suppose we could show that
$\t=Ca\Phi-D\Phi\geq 0$ on $U_\d$ with
$a=(\sup_{U} |\na\Phi| + 1)$. Since $\t(0)=0$
it would  follow that ${\pl\t\over \pl x}(0)\geq 0$
which would then yield (\ref{cl}).
\v
\noindent
To show $\t\geq 0$ on $U_\d$ it suffices to show
\v
a) $\t\geq 0$ on $\pl U_\d$

b) $\D'\t\leq 0$

\v\noindent
Now $\Phi\geq 0$ and $D\Phi=0$ if $x=Re(z_{m})=0$. Thus, to achieve a),
we try to modify $\t$ as follows: $\t= Ca\Phi+{a\over \d}|z|^2-D\Phi$. Since $D\Phi=0$
when $x=0$, we see that if $D\Phi\not=0$ on $\pl(U_\d)$, we must have
$|z|^2=\d$. Thus condition a) holds. Note
that we still have $\t(0)=0$ and we have not changed the value of
${\pl\t\over \pl x}(0)$.

\v

Next
we check b):

\v

\be
\label{one}
\D'\Phi= (g')^{\al\bar\b}( g'_{\bar\b\al}-g_{\bar\b\al})=m-\s_{\al=1}^{m}{1\over
\l_\al}
\ee
where the $\l_\al$ are the eigenvalues of $ g'$ with respect to $g$.

\v

Now  $\log\det\, g'= \log F+\log\det \,g $ implies $
(g')^{\al\bar\b}(Dg_{\bar\b\al}+\pl_\al\pl_{\bar\b}D\Phi) =  D\log F+ D\log\det g
$
so
\be
\label{two}
|\D' (D\Phi)|\ \leq \ K(1+\s{1\over \l_\al})\ \ \ {\rm and} \ \ \
\D' |z|^2\ \leq \ 2\s{1\over \l_\al}
\ee
Here $K$ is a constant depending only on $g$ (more precisely, $Dg_{\bar\b\al}$) and
on $\sup_M|\na\log F|$. Thus
$$\D'\left((K+{2\over \d})a\Phi+{a\over\d}|z|^2 -D\Phi\right)\leq Ca(m+1)$$
We are trying to show that $\D'\t\leq 0$ so this
is not quite what we wanted. Thus we modify one more time: Let
$C=K+{2\over \d}$
and define
\be
\label{mod}
\t=Ca\Phi+{a\over\d}|z|^2+Ca[x-Nx^2]-D\Phi
\ee
where $N$ is a constant, to be chosen later.
Since $x-Nx^2\geq 0$ if $\d$ is small enough, we still have $\t\geq 0$ on $\pl U_\d$.
The value ${\pl\t\over \pl x}(0)$ is replaced by
${\pl\t\over \pl x}(0)+C_3a$, which does not affect the estimate we
want. Since $\D' x^2= g^{m,\overline {m}}\geq {1\over\sup\l_\al}$ (the diagonal
entries of a hermitian matrix are greater than the smallest eigenvalue of the matrix), we have
$$
\D'\Phi+\D'(x-Nx^2)\leq m-\s{1\over
\l_\al}-{N\over
\sup_\al\l_\al}\leq m-N^{1\over m}(\det\, g')^{-1}
$$
Since $\det\, g'=F\,\det g$ we have,
 $m-N^{1\over m}(\det
\,g')^{-1}< -3m$ provided $N$ satisfies the bound: $N^{1\over m}> 2m(\sup_MF)2^{m}$.
Now we choose $\d_g$ as before, but also satisfying the condition $\left({1\over
\d}\right)^{1\over m}>2m(\sup_MF)2^{n+1}$. Thus
\be
\label{three}
\D'\Phi+\D'(x-Nx^2)\ \leq\ -3m\leq -2(m+1)\ \ \ {\rm if} \ \ \ \left({1\over
\d}\right)^{1\over n+1}> N^{1\over
m}>2m(\sup_MF)2^{m}\ .
\ee
\v
Now we obtain from (\ref{mod}):
$$
\D'\t\ = \ \D'(2Ca\Phi+{a\over\d}|z|^2+{Ca}[x-Nx^2]-D\Phi)\
$$
$$ \leq \ Ca(m+1)- 2Ca(m+1)\leq \ 0
$$
This shows $\t\geq 0$ and $\t(0)=0$. Thus
$$
{\pl\over\pl x}(D\Phi)\ \leq \ 2Ca\sup_{\pl M}|\na\Phi|+Ca\leq 2Ca(1+\sup_{\pl M}|\na\Phi|)
$$
This establishes (\ref{cl}) and the lemma is proved.

\medskip
We return to the proof of Theorem \ref{degenerateexistence}.
We rewrite  equation (\ref{MAFs}) as follows:
\bea
\label{MAFs1}
(\Omega^\e+{i\over 2}\ddb(\Phi_s+(s-1)\log\|\si\|))^m
= F_s^{(\e)}(\Omega^\e)^m,
\eea
We apply Lemma \ref{C2barrier}
to the equation (\ref{MAFs1}).
The function $|\na\log F|$ in Lemma \ref{C2barrier}
vanishes  when $F=F_s^{(\e)}$ and
we conclude that, in the flat boundary case,
$|\D\Phi_s|$ is bounded
uniformly on $\pl M$.
Thus, in this case we have
\be
|\Delta\Phi_s(z)|
\leq
C\,\|\sigma(z)\|^{-2\e A_2},
\qquad z\in M\setminus E,
\ee
with $C, A_1, A_2$ all independent of $s$ and $z$.
This completes step (d).

\v

Now choose a subsequence $\Phi_{s_j}$ converging in $C^\al$
(resp. $C^{1,\al}$) to a function $\Phi$, uniformly on
compact subsets of $ M\setminus E$.
Now equations (\ref{MA00-2}) and
(\ref{Fcondition1}) imply
\bea
\int_{M}(\Omega_s+{i\over 2}\ddb\Phi_s)^m
=
\int_{M}F_s\Omega_s^m\to 0
\quad{\rm as}\ \ s\to 0\ .
\eea
Thus the Chern-Levine-Nirenberg inequality implies
that the function $\Phi$ satisfies (\ref{MA00}) and
Theorem 2 is proved.

\medskip
Although it is not stated explicitly in
Theorem \ref{degenerateexistence},
it may be worth noting that the gradient and the $C^{1,\al}$
norms (in the case of flat boundary) of the solution $\Phi$
satisfy growth estimates near the divisor $E$
of the form $\|\sigma (z)\|^{-\e A_3}$,
where $\e A_3$ is a geometric quantity depending only on $\O_0$ and $E$.

\section{Test configurations and geodesic rays}
\setcounter{equation}{0}

We address now the problem of constructing geodesic rays associated to
a given test configuration.

\subsection{Construction of the degenerate form $\Omega_0$}

We adhere to the notation of \S 1. Thus $\r:\C^\times\ra\Aut
(\cL\ra\cX\ra\C)$ is a test configuration
for $L\ra X$,
$h_0$ is a metric on $L$ and $\o_0>0$
is its curvature.

\begin{lemma}\label{Omegazero}
There exists a metric $H_0$ on $p^*\cL^m$
with curvature $\O_0$ satisfying the following.
\begin{enumerate}
\item $\O_0\geq 0$ on $\ti\cX_D$ and
$\O_0>0$ on $\ti\cX_D^\times$.
\item $H_0$ and $\O_0$ are rotation invariant.
\item $\O_0|_{X}=\o_0$.
\end{enumerate}
\end{lemma}

\v\noindent
{\it Proof of Lemma \ref{Omegazero}.}
The Donaldson imbedding theorem
(see Lemma 7 of (\cite{PS07}))
says that for some $k>0$ there is an
imbedding

\be
I:(\cL^k\ra\cX\ra\C)\hookrightarrow
(O(1)\times\C\ra\P^N\ra\C)
\ee
Let $H_1$ be a rotation invariant  metric on $O(1)\times\C$
with positive curvature. For example,
we can take $H'= h_{FS}e^{-|w|^2}$
where $w$ is the parameter in $\C$.
Let $H_2=(I\circ p)^*H_1^{1/k}$, and
$\O_2$ its curvature.
Let $\o_2=\O_2|_X$. Since $H_2$ and $h_0$
are two metrics on the same line bundle
$L\times D^\times$ there is a smooth
function $\Psi: \ti\cX^\times = L\times D^\times\ra \R$
such that $H_2=h_0e^{-\Psi(x,w)}$.

Let $\eta:D^\times\ra [0,1]$ be a smooth
function such that $\eta(w)=1$
if $|w|\leq {1\over 3}$ and $\eta(w)=0$
if $|w|\geq {2/3}$. Let
$$ H_3=h_0e^{-\eta(w)\Psi(x,w)}
$$
Then $H_3=h_0$ if $|w|\geq {2\over 3}$ and
$H_3=H_2$ if $|w|\leq {1\over 3}$. Moreover,
the curvature of $H_3$ is positive on
the fibers $X\times\{w\}$ for all
$w\in D^\times$. Let $\al>0$ be a large
positive number and define
$H=H_3e^{-\al(|w|^2-1)}$. Although $H$
is defined as a metric on
$p^*\cL^m|_{\ti\cX^\times _D}$ it clearly
extends to a metric on ${\ti\cX^\times _D}$
which has all the desired properties.
This proves Lemma \ref{Omegazero}.

\subsection{Construction of a non-degenerate form $\Omega$}

To apply Theorem \ref{degenerateexistence} to the geodesic ray equation,
we need a divisor $E$ over $\tilde{\cal X}_D$ with $\Omega_0-\e[E]$ a K\"ahler
class. We construct here such a divisor.

\begin{lemma}
\label{divisorE}
There exists an effective divisor $E$ supported on the central
fiber of $\ti\cX_D$ and an $\e>0$ with the following
property: $\O-\e[E]$ is a \K class
(where $[E]$ is the current of integration).
\end{lemma}

\noindent More precisely, we claim that there is a smooth metric $H$
on the line bundle $O(E)$ and a positive \K metric $\O_\e$ on
$\ti\cX$, such that

\be\label{assumption1}
\O_0-\e[E]\ = \ \O_\e-\e{i\over 2} \ddb\log \|\sigma\|^2
\ee
where $\sigma$ is the canonical section
of $O(E)$ which vanishes on $E$, and
$\|\sigma(z)\|^2=|\sigma(z)|^2H(z)$.
\v\noindent

{\it Proof.} In Lemma 2 of \cite{reg} it was shown that there
is a line bundle $\cM\ra \ti\cX$ and an integer $k>0$
such that

\v

1) $p^*\cL^{mk}\otimes\cM\ra\ti\cX$ is ample, that it, it has a metric $h_k$
of positive curvature.

2) There is a meromorphic section $\m:\ti\cX\ra\cM$ whose retriction
to $\ti\cX^\times$ is holomorphic

\hskip .19in and non-vanishing.
\v

In fact the proof shows that
$\m={1\over s}$, where $s$
is a global section of $\cM^{-1}$ which is non-vanishing outside the
central fiber (see the paragraph before (3.8) in \cite{reg}).

\v
Now let $\O_k>0$ be the curvature of
of the metric $h_k$ and let $E=\{s=0\}$,
so that $\cM=O(-E)$. Let $H=h_k^{-1}H_0^k$. Then

$$ -{i\over 2}\ddb\log h\ = \ k\O_0\ - \ \O_k
$$
If we let $\e={1\over k}$ and $\O_\e={1\over k}\O_k$ and $s=\sigma$, then
we obtain (\ref{assumption1}).

\v\v

\subsection{Existence of geodesic rays}
We can now prove Theorem 3:
Consider the Dirichlet problem  (\ref{MA00}) for the
manifold
 $M=\tilde{\cal X}_D$ where
$\Omega_0$ is the non-negative form constructed in Lemma \ref{Omegazero},
and $\varphi=0$.
Letting the divisor $E$
be the divisor constructed in Lemma \ref{divisorE}, we readily see
that all the hypotheses of Theorem
\ref{degenerateexistence}, part (a) are satisfied.
Moreover, the manifold $M$ has a locally flat boundary, so
Theorem \ref{geodesicray} follows from Theorem \ref{degenerateexistence}.

\bigskip
\noindent
{\bf Acknowledgements}

The authors would particularly like to thank Jian Song for
calling their attention to the work of Tsuji in this context and for some
helpful conversations. They would like to thank Zbigniew Blocki
for communicating to them his preprint \cite{B09}
prior to publication, and Steve Zelditch for some stimulating discussions.
They would also like to thank Ovidiu Munteanu for providing them
with the references \cite{CY86} and \cite{K83}.

\noindent Department of Mathematics \\
Columbia University, New York, NY 10027

\medskip

\noindent Department of Mathematics \\
Rutgers University, Newark, NJ 07102

\end{document}